\documentclass[10pt, twoside]{amsart}
\advance\oddsidemargin by -1.0cm
\advance\evensidemargin by -1.0cm
\advance\topmargin by -1.0cm 
\usepackage{amssymb}
\usepackage{amstext}
\usepackage{amscd}   

\theoremstyle{plain}
\newtheorem{cor}{Corollary}[section]
\newtheorem{lem}{Lemma}[section]
\newtheorem{thm}{Theorem}[section]
\newtheorem{prop}{Proposition}[section]
\theoremstyle{definition}

\newcommand{\R}{\ensuremath{\mathbb{R}}}

\newcommand{\tr}{\ensuremath{\mathrm{tr}\,}}

\newcommand{\Ric}{\ensuremath{\mathrm{Ric}}}

\newcommand{\be}{\begin{equation}}
\newcommand{\ee}{\end{equation}}
\newcommand{\bea}{\begin{eqnarray}}
\newcommand{\eea}{\end{eqnarray}}
\newcommand{\bean}{\begin{eqnarray*}}
\newcommand{\eean}{\end{eqnarray*}}

\begin{document}

\title[Integrability conditions for almost Hermitian and almost 
K\"ahler $4$-manifolds]{Integrability conditions for almost Hermitian and almost K\"ahler $4$-manifolds}

\address{
{\normalfont\ttfamily 
Institute of Mathematics \newline
Humboldt-Universit\"at zu Berlin\newline
Office: Rudower Chaussee 25\newline
D-10099 Berlin, Germany\newline
kirchber@mathematik.hu-berlin.de}}
\author[K.-D. Kirchberg]{K.-D. Kirchberg}
\maketitle 

\begin{center}
(Version of \today)
\end{center}

\begin{abstract}
\noindent If
 $W_+$ denotes the self dual part of the Weyl tensor of any K\"ahler
4-manifold
and $S$ its scalar curvature, then the relation $|W_+|^2 = S^2/6$ is
 well-known. For any almost K\"ahler 4-manifold with $S \ge 0$, this condition
forces the K\"ahler property. A compact almost K\"ahler 4-manifold is already
K\"ahler if it satisfies the conditions $| W_+ |^2 = S^2/6$ and $\delta W_+=0$ and also
if it is Einstein and $| W_+|$ is constant. Some further results of this type are proved.
An almost Hermitian 4-manifold $(M,g,J)$ with $\mathrm{supp} (W_+)=M$ is already K\"ahler
if it satisfies the condition $| W_+ |^2 = 3 (S_{\star} - S/3)^2 /8$
together with $|\nabla W_+ | = | \nabla |W_+||$ or with $\delta W_+ + \nabla \log | W_+ |
\lrcorner W_+ =0$, respectively. The almost complex structure $J$ enters here explicitely via the 
star scalar curvature $S_{\star}$ only.\\

\thispagestyle{empty}

\noindent 2002 Mathematics Subject Classification: 53B20, 53C25
\end{abstract}
\bigskip

\setcounter{section}{-1}
\noindent\section{\bf Introduction}

This paper is related to the following general question: Which 
curvature properties of an almost Hermitian (almost K\"ahler) manifold
$(M,g,J)$ imply that $(M,g,J)$ is in fact a Hermitian manifold
or a K\"ahler manifold, respectively? On can also consider the stronger
problem which curvature properties of the corresponding oriented Riemannian
manifold $(M,g)$ force that $(M,g,J)$ is already Hermitian or K\"ahler,
respectively. By definition, curvature properties of this kind imply
the integrability of the almost complex structure $J$ and, hence, 
are called integrability conditions for almost Hermitian (almost K\"ahler)
manifolds. Here we investigate these problems in dimension $4$ and we formulate
integrability conditions where the Weyl tensor enters essentially. This was 
motivated by the fact that the very special form of the self-dual part
of the Weyl tensor of a K\"ahler $4$-manifold yields several necessary 
integrability conditions. Hence, it is very natural to ask which of these 
conditions are also sufficient. In more detail, 
for any K\"ahler 4-manifold $(M, g,J)$ with scalar curvature $S$ and self-dual
part $W_+ : \Lambda^2_+ \to \Lambda^2_+$ of the Weyl tensor $W$, the relations
\begin{equation} \label{gl-neu1}
|W_+ |^2 = \frac{S^2}{6}  \ , 
\end{equation}

\begin{equation} \label{gl-neu2}
(\det (W_+))^2 = \frac{1}{54} |W_+|^6  \ ,
\end{equation}

\begin{equation} \label{gl-neu3}
| \nabla W_+ |^2 = \frac{1}{6} | \nabla S|^2
\end{equation}

are well-known. From (\ref{gl-neu1}) and (\ref{gl-neu3}) we obtain the equation
\begin{equation} \label{gl-neu4}
|\nabla W_+ | = | \nabla | W_+ \|
\end{equation} 

\renewcommand{\thefootnote}{\fnsymbol{footnote}}
which is valid on the open subset $M_+ \subset M$ where $W_+$ does not vanish.
 Moreover, one has the equation
\begin{equation} \label{gl-neu5}
\delta W_+ + \nabla \log | W_+ | \lrcorner W_+ =0
\end{equation}

on $M_+$. Thus, the K\"ahler property imposes strong conditions on $W_+$.  This
paper deals with the converse question which curvature conditions already
imply the K\"ahler property of an almost Hermitian 4-manifold $(M,g,J)$ or,
more specially, of an almost K\"ahler 4-manifold.\\

In the almost Hermitian case, we replace (\ref{gl-neu1}) by the basic condition
\begin{equation} \label{gl-neu6}
 | W_+ |^2 = \frac{3}{8} (S_{\star} - \frac{S}{3} )^2
\end{equation}

containing the star scalar curvature $S_{\star}$ and being equivalent to (\ref{gl-neu1}) 
if $S_{\star} = S$. If $M_+$ is dense in $M$ \ $(\mathrm{supp} (W_+)=M)$, we
show that (\ref{gl-neu4}), (\ref{gl-neu6}) and also (\ref{gl-neu5}), (\ref{gl-neu6}) are sufficient conditions for 
an almost Hermitian 4-manifold to be K\"ahler (Theorem \ref{thm-3-1} and 
Theorem \ref{thm-3-3}). Combining the condition
\begin{equation} \label{gl-neu7}
\delta W_+ =0
\end{equation}

with (\ref{gl-neu6}) we obtain that $J$ is integrable on $M_+$. Thus, every almost
Hermitian 4-manifold $(M,g,J)$ with (\ref{gl-neu6}), (\ref{gl-neu7}) and $\mathrm{supp} (W_+)=M$ is
already Hermitian (Theorem \ref{thm-3-2}). In this case, the restriction of $g$ 
to $M_+$ is conformally equivalent to a metric $\bar{g}$ on $M_+$ such that 
$(M_+ , \bar{g}, J)$ is a K\"ahler manifold. 
This result is related to a theorem which was published
by A.~Derdzinski already in 1983 \cite{8}. He considered oriented Riemannian
4-manifolds for which $W_+$ has at most two eigenvalues at every point
and satisfies (\ref{gl-neu7}). Then he proved that on $M_+$ the metric $g$ is
conformally equivalent to a metric $\bar{g}$ which is locally K\"ahler (\cite{6}, 
16.67). The relation to this result follows from the fact that (\ref{gl-neu6}) implies 
that $W_+$ has exactly two eigenvalues on $M_+$ or, equivalently, satisfies 
the conformally invariant equation (2). A result similar to the mentioned
result of Derdzinski is Theorem \ref{thm-3-4} which states that an oriented
Riemannian 4-manifold $(M,g)$ with $W_+ \not= 0$ everywhere is locally K\"ahler
with respect to local complex structures that are compatible with the given 
orientation if and only if $W_+$ satisfies the conditions (\ref{gl-neu2}) and (\ref{gl-neu4})
or (\ref{gl-neu2}) and (\ref{gl-neu5}), respectively. In Section 4, we finally consider integrability
conditions  for almost K\"ahler 4-manifolds. In connection with the conjecture
of S.I.~Goldberg \cite{11}, interesting results in this direction were proved already
(\cite{2}-\cite{4}, \cite{9}, \cite{10}, \cite{12}, \cite{14}-\cite{18}). The Goldberg conjecture states
that any compact almost K\"ahler Einstein manifold is in fact K\"ahler. For non-negative scalar curvature $S$, it was shown by K.~Sekigawa \cite{18} 
that this conjecture is true. For $S <0$, no proof is known so far. There are
attempts to construct counterexamples of this part of the Goldberg conjecture
(\cite{1}, \cite{13}). A result which also uses the assumption $S \ge 0$ is
our Theorem \ref{thm-4-1}. It states that an almost K\"ahler 4-manifold with $S
\ge 0$ is K\"ahler if the curvature  condition (\ref{gl-neu1}) is fulfilled. 
In the compact case, our basic result is Theorem \ref{thm-4-2}. By this theorem,
a compact almost K\"ahler 4-manifold $(M,g,J)$ is K\"ahler if and only if
it satisfies condition (1) together with $Q(J)=0$, where the number $Q(J)$
is an obstruction to the K\"ahler property for any compact almost Hermitian
manifold. $Q(J)$ vanishes, for example, if the Weyl tensor $W$ is harmonic
 $(\delta W=0)$. An application of Theorem \ref{thm-4-2} is Theorem \ref{thm-4-3},
which involves two assertions. The first one states that a compact almost
K\"ahler 4-manifold is K\"ahler if and only if it satisfies (\ref{gl-neu1}) together
with
\begin{equation} \label{gl-neu8}
\delta W_+ + \Theta = 0 \ ,
\end{equation}

where $\Theta$ is a certain tensor field of type $(2,1)$ depending on $J$ purely
algebraically. By the second assertion of Theorem \ref{thm-4-3}, a compact almost
K\"ahler 4-manifold is a K\"ahler manifold of constant scalar curvature if and
only if it has the properties (\ref{gl-neu1}) and (\ref{gl-neu7}). We remark that condition
(\ref{gl-neu7}) and also condition (\ref{gl-neu8}) generalize the Einstein condition. Thus, our
Theorem \ref{thm-4-3} is related to the Goldberg conjecture. For example, it
yields an improvement of Corollary F in \cite{17} which 
asserts that any compact almost K\"ahler Einstein 4-manifold of negative scalar
curvature is K\"ahler if it satisfies (\ref{gl-neu1}). The final theorems are also related to the 
Goldberg conjecture. We show that a compact almost K\"ahler Einstein 4-manifold
is already K\"ahler if $| W_+|$ is constant (Theorem \ref{thm-4-4}). Moreover,
by Theorem \ref{thm-4-5}, a compact 
almost K\"ahler 4-manifold with constant negative scalar curvature $S$ is
K\"ahler  if it satisfies equation (\ref{gl-neu7}) together with the necessary condition
\begin{equation} \label{gl-neu9}
\mathrm{det} (W_+) = \frac{S^3}{108} \ . 
\end{equation}

In particular, this implies that a compact almost K\"ahler Einstein 4-manifold with $S <0$ and 
(\ref{gl-neu9}) is in fact K\"ahler (Corollary \ref{cor-4-4}).

\section{Preliminaries}

Let $(M,g)$ be any Riemannian $n$-manifold. Then, for any endomorphism of the 
tangent bundle $A:TM \to TM$, we denote by $A^* :TM \to TM$ the corresponding
adjoint endomorphism defined by the characteristic property $g(AX,Y)=g(X, A^* Y)$.
The almost Hermitian structure of an almost Hermitian $n$-manifold $(M,g,J)$ is
characterized by the basic properties
\begin{equation*} \label{gl-01}
J^2 = - 1 \quad , \quad J^* = - J
\end{equation*}

of the almost complex structure $J$ which imply that
\begin{equation*} \label{gl-02}
g(JX, JY) = g(X,Y)
\end{equation*}

for all vector fields $X,Y$. The fundamental 2-form $\Omega$ is defined by
$\Omega (X,Y):=g(JX,Y)$ and the canonical orientation is given by the volume
$n$-form $\omega := \frac{1}{m!} \Omega^m \ (n=2m)$. As usual, by $\nabla$ we denote
the Levi-Civita covariant derivative corresponding to $g$. We use the notation
\begin{displaymath}
\nabla^2_{X,Y} := \nabla_X \circ \nabla_Y - \nabla_{\nabla_X Y} 
\end{displaymath}

for the tensorial covariant derivatives of second order such that the Riemannian
curvature tensor $R$ is given by $R(X,Y)Z= \nabla^2_{X,Y} Z - \nabla^2_{Y,X} Z$.
In our paper, the Ricci tensor is the endomorphism $\Ric :TM \to TM$ locally 
defined by $\Ric (X) := R(X, X_k)X^k$. Here $(X_1, \ldots , X_n)$ denotes a local 
frame of vector fields and $(X^1 , \ldots, X^n)$ the associated frame defined by
$X^k := g^{kl} \cdot X_l \ (k=1, \ldots , n)$, where $(g^{kl})$ is the inverse 
of the matrix $(g_{kl})$ with $g_{kl} := g(X_k, X_l)$. In case of an orthonormal
frame, we have $X^k = X_k \ (k=1, \ldots , n)$ then. For any almost Hermitian
manifold $(M, g,J)$, the star Ricci tensor $\Ric_{\star} : TM \to TM$ is defined
by $\Ric_{\star} (X):= R(JX, JX_k) X^k$. The first Bianchi identity yields the 
equation
\begin{equation} \label{gl-03}
\Ric_{\star} = - \frac{1}{2} R(J X_k , X^k) \circ J
\end{equation}

which implies
\begin{equation} \label{gl-04}
(\Ric_{\star})^* = - J \circ \Ric_{\star} \circ J \ . 
\end{equation}

We also use the notations
\begin{eqnarray*}
\Ric^{\pm} &:= &\frac{1}{2} (\Ric \mp J \circ \Ric \circ J) \ , \\
\Ric^{\pm}_{\star} & := & \frac{1}{2} (\Ric_{\star} \pm (\Ric_{\star})^*) \ . 
\end{eqnarray*}

By definition, then we have
\begin{equation} \label{gl-05} 
\Ric = \Ric^+ + \Ric^- \quad , \quad \Ric_{\star} = \Ric_{\star}^+ +
\Ric^-_{\star} \ , 
\end{equation}
\begin{equation} \label{gl-06} 
(\Ric^{\pm})^* = \Ric^{\pm} \quad , \quad (\Ric^{\pm}_{\star} )^* = \pm
\Ric_{\star}^{\pm} \ . 
\end{equation}

Moreover, the endomorphisms $\Ric^+$ and $\Ric^+_{\star}$ ($\Ric^-$ and
$\Ric_{\star}^-$) commute (anticommute) with $J$, i.e., it holds that
\begin{equation} \label{gl-07} 
[\Ric^+ , J] =[ \Ric_{\star}^+ , J] = \{ \Ric^- , J \} = \{ \Ric_{\star}^- , 
J \} =0 \ , 
\end{equation}

where $[A,B] := A \circ B - B \circ A$ denotes the commutator and 
$\{ A, B \} := A \circ B + B \circ A$ the anticommutator of endomorphisms
$A,B$. The Ricci form $\rho$ and the star Ricci forms $\rho_{\star}, 
\rho^+_{\star} , \rho^-_{\star}$ are defined by
\begin{displaymath}
\rho (X,Y)  := g(( \Ric^+ \circ J) X,Y) \quad , \quad \rho_{\star} (X,Y) :=
g(( \Ric_{\star} \circ J)X, Y) \ , 
\end{displaymath}
\begin{displaymath}
\rho^{\pm}_{\star} (X,Y) := g(( \Ric^{\pm}_{\star} \circ J) X,Y) \ . 
\end{displaymath}

By definition, we have the decomposition
\begin{equation*} \label{gl-08} 
\rho_{\star} = \rho^+_{\star} + \rho^-_{\star} \ . 
\end{equation*}

Moreover, (\ref{gl-07}) yields
\begin{equation*} \label{gl-09} 
\rho (JX , JY)= \rho (X,Y) \quad , \quad \rho^{\pm}_{\star} (JX , JY) =
\pm \rho^{\pm}_{\star} (X,Y) \ . 
\end{equation*}

Besides the scalar curvature $S := \tr (\Ric) = \tr (\Ric^+)$ we
also consider the star scalar curvature $S_{\star} := \tr (\Ric_{\star})=
\tr (\Ric_{\star}^+)$. The tensor field $\tilde{R}$ defined by
\begin{displaymath}
\tilde{R} (X,Y) := \frac{1}{4} [R (X,Y) - R(JX, JY) , J] \circ J
\end{displaymath}

has the properties
\begin{equation*} \label{gl-10} 
\tilde{R} (X,Y)^* = - \tilde{R} (X,Y) = \tilde{R} (Y, X) \ , 
\end{equation*}
\begin{equation*} \label{gl-11} 
\tilde{R} (JX, JY) = - \tilde{R} (X,Y) \ , 
\end{equation*}
\begin{equation*} \label{gl-12} 
\{ \tilde{R} (X,Y) , J \} = 0 \ . 
\end{equation*}

Furthermore, we use the decomposition
\begin{equation*} \label{gl-13} 
\tilde{R} (X,Y) = \tilde{R}^+ (X, Y) + \tilde{R}^- (X, Y)
\end{equation*}

with $\tilde{R}^{\pm} (X,Y) := \frac{1}{2} (\tilde{R} (X,Y) \pm 
\tilde{R} (JX, Y) \circ J)$. Then
$\tilde{R}^+$ and $\tilde{R^-}$  have the additional properties
\begin{equation*} \label{gl-14} 
\tilde{R}^{\pm} (JX, Y) \circ J= \pm \tilde{R}^{\pm} (X,Y) \ . 
\end{equation*}

It is well-known that $\tilde{R}^-$ is already determined by the Weyl tensor $W$, 
i.e., we have
\begin{equation*} \label{gl-15} 
\tilde{R}^- (X,Y) = \tilde{W}^- (X,Y) \ . 
\end{equation*}

\newcommand{\wric}{\widetilde{\Ric}}
We introduce the curvature endomorphism $\wric$ defined by $\wric (X) :=
\tilde{R} (X, X_k)X^k$ and the function $\tilde{S} := \tr (\wric)$. A
direct calculation yields the identity
\begin{equation} \label{gl-16} 
\wric = \frac{1}{2} (\Ric_{\star}^+ - \Ric^+)
\end{equation} 

implying
\begin{equation} \label{gl-17} 
\tilde{S} = \frac{1}{2} (S_{\star} - S) \ , 
\end{equation}
\begin{equation*} \label{gl-18} 
(\wric)^* = \wric \quad , \quad [ \wric , J] = 0 \ . 
\end{equation*}

For any skew symmetric endomorphism $A:TM \to TM \ (A^* = - A)$ we define skew 
symmetric endomorphisms $R(A), \tilde{R} (A), W(A)$ by $R(A) := R(AX_k, X^k), 
\tilde{R} (A) := \tilde{R} (AX_k , X^k)$ and $W(A):= W(AX_k, X^k)$. Then,
the equations
\begin{equation} \label{gl-19} 
R(A) =W(A) + \frac{2}{n-2} (\{ \Ric , A\} - \frac{S}{n-1} A) \ , 
\end{equation}
\begin{equation} \label{gl-20} 
\tilde{R} (A) = \frac{1}{4} [R(A) + R(J \circ A \circ J) , J] \circ J \ , 
\end{equation}
\begin{equation*} \label{gl-21} 
\{ \tilde{R} (A) , J \} = 0 \ , 
\end{equation*}
\begin{equation*} \label{gl-22} 
\tilde{R} (J \circ A \circ J) = \tilde{R} (A)
\end{equation*}

are valid for all skew symmetric endomorphisms $A$. We remark, that usually 
$R, \tilde{R}$ and $W$ are considered as endomorphisms of the bundle 
$\Lambda^2 := \Lambda^2 T^* M$. But here we use the canonical isomorphism
between skew symmetric endomorphisms and 2-forms given by
\begin{equation} \label{gl-23} 
\Omega_A (X,Y) = g (AX, Y) \ , 
\end{equation}

where $\Omega_A$ denotes the 2-form corresponding to the skew symmetric
endomorphism $A$. For example, we have  $\Omega = \Omega_J$ then and $W
(\Omega_A)$ and $W(A)$ are related by
\begin{equation*} \label{gl-24} 
W(\Omega_A)(X,Y) = g(W(A)X, Y) \ . 
\end{equation*}

Using (\ref{gl-03}) and (\ref{gl-19}) we obtain
\begin{equation} \label{gl-25} 
W(J) = \frac{2S}{(n-1)(n-2)} J + 2 (\Ric_{\star} - \Ric^+) \circ J +
\frac{n-4}{n-2} \{ \Ric , J\} \ . 
\end{equation}

In the following, for any endomorphisms $A,B :TM \to TM$, 
we use the scalar product $\langle A,B \rangle$ defined by
\begin{displaymath}
\langle A,B \rangle = \frac{1}{n} \tr (A^* \circ B)= \frac{1}{n} g(AX_k, BX^k)
\end{displaymath}

which differs from the usual one by the dimension factor $1/n$. According to
 (\ref{gl-23}) we define the scalar product of 2-forms such that
\begin{equation*} \label{gl-26} 
\langle \Omega_A , \Omega_B \rangle = \langle A, B \rangle \ . 
\end{equation*}

Now, let us consider the almost K\"ahler case. An almost K\"ahler manifold
is an almost Hermitian manifold $(M,g,J)$ with closed fundamental 2-form
$\Omega = \Omega_J$, i.e., with
\begin{equation*} \label{gl-27} 
d \Omega =0 \ . 
\end{equation*}

It is well-known, that
\begin{equation*} \label{gl-28} 
\delta \Omega =0
\end{equation*}

follows, i.e., $\Omega$ is also co-closed then. Furthermore, the almost complex structure $J$ of any almost K\"ahler manifold satisfies the so-called
quasi K\"ahler condition
\begin{equation*} \label{gl-29} 
\nabla_{JX} J= \nabla_X J \circ J \ . 
\end{equation*}

In the almost K\"ahler case, the tensor $\wric$ has the special form (\cite{12}, eq. (41))
\begin{equation*} \label{gl-30} 
\wric = - \frac{1}{4} \nabla_{X_k} J \circ \nabla_{X^k} J \ . 
\end{equation*}

This implies $\wric \ge 0$ and, moreover, 
\begin{equation} \label{gl-31} 
\tilde{S} = \frac{n}{4} | \nabla J |^2 \ , 
\end{equation}

where, according to the definition of the scalar product above, the function 
$| \nabla J |^2$ is locally given by
\begin{displaymath}
| \nabla J |^2 = \langle \nabla_{X_k} J , \nabla_{X^k} J \rangle = - \frac{1}{n}
\tr (\nabla_{X_k} J \circ \nabla_{X^k} J) \ . 
\end{displaymath}

By (\ref{gl-17}) and (\ref{gl-31}), we have
\begin{equation} \label{gl-32} 
S_{\star} - S = \frac{n}{2} | \nabla J |^2
\end{equation}

in the almost K\"ahler case.

\section{The Weyl tensor of an almost Hermitian 4-manifold.}

For any oriented Riemannian 4-manifold $(M,g)$, the Hodge operator $*$ acts
as an involution on the bundle $\Lambda^2 := \Lambda^2 T^* M$. This yields
the orthogonal splitting
\begin{equation*} \label{gl-33} 
\Lambda^2 = \Lambda^2_+ \oplus \Lambda^2_-
\end{equation*}

with rank $(\Lambda^2_{\pm})=3$, where $\Lambda^2_+ (\Lambda^2_-)$ denotes the 
eigen-subbundle of this involution to the eigenvalue $1 \ (-1)$. We use the
notations $P_+$ and $P_-$ for the corresponding bundle projections
$(P_{\pm} (\Lambda^2) = \Lambda^2_{\pm})$. Since the Weyl tensor $W$ considered
as an endomorphism of $\Lambda^2$ commutes with $P_{\pm}, W$ decomposes in the 
form $W= W_+ \oplus W_-$, where $W_+ : \Lambda^2_+ \to \Lambda^2_+ (W_- : 
\Lambda^2_- \to \Lambda^2_-)$ is called the self-dual (anti-self-dual) part
of $W$. $W_+$ and $W_-$ are self-adjoint and traceless endomorphisms. We
summarize some well-known basic facts concerning this splitting in the 
following lemma (see \cite{5}).

\begin{lem} \label{lem-2-1} \mbox{} \\
(i) For any section $\Omega_A$ of $\Lambda^2_+$ and any section $\Omega_B$ of
$\Lambda^2_-$, $A$ and $B$ commute $([A,B]=0)$.\\
(ii) If $\Omega_J \in \Gamma (\Lambda^2_+) \ (\Omega_J \in \Gamma (\Lambda^2_-))$
is any unit section $(|\Omega_J |=1)$, then $(M,g,J)$ is an almost Hermitian
manifold and we have $\frac{1}{2} \Omega_J \wedge \Omega_J = \omega$ $(\frac{1}{2}
\Omega_J \wedge \Omega_J = - \omega)$, where $\omega$ is the volume $4$-form 
defined by the metric $g$ and the given orientation.\\
(iii) $\Lambda^2_+$ and $\Lambda^2_-$ are endowed with a canonical orientation, 
where a local orthonormal frame $(\Omega_I, \Omega_J, \Omega_K)$ of
$\Lambda^2_+$ or $\Lambda^2_-$, respectively, is positively oriented if and 
only if the corresponding frame $(I,J,K)$ satisfies the quaternionic relations
$I^2 = J^2 = K^2 = - 1$, $I \circ J \circ K = - 1$.
\end{lem}

We now consider an almost Hermitian 4-manifold $(M,g,J)$ endowed with the 
orientation given by $J$. Then the corresponding fundamental 2-form $\Omega =
\Omega_J$ satisfies the equation $* \Omega = \Omega$ implying $\Omega
\in \Gamma (\Lambda^2_+)$. Thus, we obtain the orthogonal splitting
\begin{equation} \label{gl-34} 
\Lambda^2_+ = \R \cdot \Omega \oplus P_+ (\Omega^{\perp}) \ . 
\end{equation}

Let $x \in M$ be any point and let $(X_1 , \ldots, X_4)$ be any orthonormal
frame of vector fields in a neighborhood $U$ of $x$ with the property
\begin{equation*} \label{gl-35} 
JX_1 = X_2 \quad , \quad JX_2 = - X_1 \quad , \quad JX_3 = X_4 \quad , \quad
JX_4 = - X_3 \ . 
\end{equation*}

Such a frame is called a $J$-frame. We consider the endomorphisms $I,K$ on $U$
 defined by
\begin{equation*} \label{gl-36} 
\begin{array}{l}
IX_1 = X_3 \quad , \quad IX_2 = - X_4 \quad , \quad IX_3 = - X_1 \quad , 
\quad IX_4 = X_2 \ , \\[0.5em]
KX_1 = - X_4 \quad , \quad KX_2 = - X_3 \quad , \quad KX_3 = X_2 \quad
, \quad KX_4 = X_1 \ . 
\end{array}
\end{equation*}

Then, by definition, we have
\begin{equation*} \label{gl-37} 
I^* = - I \quad , \quad K^* = - K
\end{equation*}

and the quaternionic relations
\begin{equation*} \label{gl-38} 
I^2 = J^2 = K^2 = - 1 \quad , \quad I \circ J \circ K = - 1
\end{equation*}

are satisfied. In the following, such a pair of
endomorphisms $(I,K)$ is called a local quaternionic supplement
of $J$. Moreover, it holds that 
\begin{equation*} \label{gl-39} 
|I|^2 = |J|^2 = |K|^2 =1 \quad , \quad \langle I,J \rangle  = \langle 
I,K \rangle = \langle J,K \rangle =0 \ , 
\end{equation*}
\begin{equation*} \label{gl-40} 
* \Omega_I = \Omega_I \quad , \quad * \Omega_K = \Omega_K \ . 
\end{equation*}

Thus, every triple $(\Omega, \Omega_I , \Omega_K)$ is a (negatively oriented)
local orthonormal frame of $\Lambda^2_+$, where $(\Omega_I , \Omega_K)$ is a
frame of $P_+ (\Omega^{\perp})$. Furthermore, two quaternionic supplements
$(I, K)$ and  $(I', K')$ of $J$ on a sufficiently small open neighborhood $U$ are 
related by
\begin{equation*} \label{gl-41} 
I' = \cos \alpha \cdot I - \sin\alpha \cdot K \quad , \quad K' = \sin \alpha \cdot
I + \cos \alpha \cdot K
\end{equation*}

with a function $\alpha$ on $U$.\\

For any local quaternionic supplement $(I,K)$ of $J$, we introduce the local 
endomorphisms $\Ric_{\bigtriangleup}$ and $\Ric_{\Box}$ defined by
\begin{displaymath}
\Ric_{\bigtriangleup} (X) := \Ric (IX , IX_k) X^k \quad , \quad
\Ric_{\Box} (X) := \Ric (KX, KX_k) X^k \ . 
\end{displaymath}

By definition, $\Ric_{\bigtriangleup}$ and $\Ric_{\Box}$ are the star Ricci tensors
corresponding to the local almost Hermitian structures $(g,I)$ and $(g,K)$, 
respectively.

\begin{lem} \label{lem-2-2} 
For every local quaternionic supplement $(I,K)$ of $J$, we have the relation
\begin{equation} \label{gl-42} 
\Ric_{\bigtriangleup} + \Ric_{\Box} + \Ric_{\star} = \Ric \ . 
\end{equation}
\end{lem}

\begin{proof}
Let $(X_1 , \ldots, X_4)$ be the $J$-frame according to which $(I,K)$ is defined.
Using the first Bianchi identity we calculate
\begin{eqnarray*}
\Ric_{\bigtriangleup} (X_1) &=& R(IX_1 , IX_2) X_2 + R(IX_1, IX_3)X_3 +
R(IX_1, IX_4)X_4 =\\
&=& R(X_3, - X_4)X_2 + R(X_3 , - X_1) X_3 + R(X_3, X_2)X_4 =\\
&=& R(X_1 , X_3)X_3 + R(X_4, X_2) X_3 \ . 
\end{eqnarray*}

In this way, we obtain the images $\Ric_{\bigtriangleup} (X_k), \Ric_{\Box} 
(X_k)$ and $\Ric_{\star} (X_k)$ for $k=1, \ldots, 4$. Then, we find
\begin{displaymath}
\Ric_{\bigtriangleup} (X_k) + \Ric_{\Box} (X_k) + \Ric_{\star} (X_k) =
\Ric (X_k) \quad \quad (k= 1, \ldots, 4) \ . 
\end{displaymath}
\end{proof}

According to (\ref{gl-05}) we have the decompositions
\begin{equation*} \label{gl-43}
\Ric_{\bigtriangleup} = \Ric_{\bigtriangleup}^+ + \Ric^-_{\bigtriangleup} \quad , \quad \Ric_{\Box} = 
\Ric^+_{\Box} + \Ric^-_{\Box}
\end{equation*}

with the properties corresponding to (\ref{gl-06}) and (\ref{gl-07}), respectively. Hence, 
Taking the symmetric and the skew symmetric part of (\ref{gl-42}) we see that
(\ref{gl-42}) is equivalent to the two equations
\begin{equation} \label{gl-44} 
\Ric^+_{\bigtriangleup} + \Ric^+_{\Box} + \Ric^+_{\star} = \Ric \ , 
\end{equation}
\begin{equation} \label{gl-45} 
\Ric^-_{\bigtriangleup} + \Ric^-_{\Box} + \Ric^-_{\star} = 0 \ . 
\end{equation}

Moreover, with $S_{\bigtriangleup} := \tr (\Ric_{\bigtriangleup})$ and $S_{\Box} 
:= \tr (\Ric_{\Box})$ from (\ref{gl-42}) we obtain
\begin{equation} \label{gl-46} 
S_{\bigtriangleup} + S_{\Box} + S_{\star} = S \ . 
\end{equation}

Since, the endomorphism 
$\wric$ is a multiple of the identity for any almost Hermitian 4-manifold (\cite{17}, (2.3.1)), from
(\ref{gl-16}), (\ref{gl-17}) we obtain here
\begin{equation*} \label{gl-47} 
\wric = \frac{\tilde{S}}{4} = \frac{1}{8} (S_{\star} - S) = \frac{1}{2} (\Ric_{\star}^+
- \Ric^+ ) \ . 
\end{equation*} 

For $n=4$, (\ref{gl-25}) can thus be written in the form
\begin{equation} \label{gl-48} 
W(J) = \frac{1}{2} (S_{\star} - \frac{S}{3}) J + 2 \Ric^-_{\star} \circ
J \ . 
\end{equation}

Hence, for any local quaternionic supplement $(I, K)$ of $J$, we have quite analogously
\begin{eqnarray}
W(I) &=& \frac{1}{2} (S_{\bigtriangleup} - \frac{S}{3})I + 2 \Ric^-_{\bigtriangleup}
\circ I , \label{gl-49}\\[0.5em]
W(K) &=& \frac{1}{2} (S_{\Box} - \frac{S}{3})K + 2 \Ric^-_{\Box}
\circ K , \label{gl-50}
\end{eqnarray}

We remark that, using the correspondence (19) between skew symmetric 
endomorphisms and $2$-forms, the equations
(\ref{gl-48})-(\ref{gl-50}) can equivalently be stated 
\begin{eqnarray}
W_+ (\Omega) &=& \frac{1}{2} (S_{\star} - \frac{S}{3}) \Omega + 2 \rho^-_{\star},\label{gl-51}  \\[0.5em]
W_+ (\Omega_I) &=& \frac{1}{2} (S_{\bigtriangleup} - \frac{S}{3}) \Omega_I + 2 \rho^-_{\bigtriangleup},\label{gl-52}  \\[0.5em]
W_+ (\Omega_K) &=& \frac{1}{2} (S_{\Box} - \frac{S}{3}) \Omega_K + 2 \rho^-_{\Box} \ . \label{gl-53} 
\end{eqnarray}

Using (\ref{gl-46}) the equations (\ref{gl-51})-(\ref{gl-53}) yield
\begin{equation} \label{gl-54} 
|W_+ |^2 = \frac{1}{4} (S_{\star}^2 + S_{\bigtriangleup}^2 + S_{\Box}^2 - 
\frac{S^2}{3}) + 4(| \Ric^-_{\star} |^2 + |\Ric^-_{\bigtriangleup} |^2 + | \Ric^-_{\Box} |^2) \ . 
\end{equation}

Applying (\ref{gl-03}) in the case where $J$ is replaced by $I$ or $K$, 
respectively, (\ref{gl-20}) yields the equations
\begin{eqnarray} 
\tilde{R} (I) &=& - \Ric_{\bigtriangleup} \circ I - J \circ \Ric_{\bigtriangleup}
\circ K , \label{gl-55} \\[0.5em]
\tilde{R} (K) &=& - \Ric_{\Box} \circ K - J \circ \Ric_{\Box}
\circ I \ .  \label{gl-56} 
\end{eqnarray}

This shows that $\{ \tilde{R} (I) , J \} = \{ \tilde{R} (K) , J \} =0$. 
Hence, $\tilde{R} (I)$ and $\tilde{R} (K)$ correspond to local sections
of $P_+ (\Omega^{\perp})$ implying that these endomorphisms are linear
combinations of $I$ and $K$. Thus, using (\ref{gl-45}) from (\ref{gl-55}), 
(\ref{gl-56}) we obtain
\begin{eqnarray}
\tilde{R} (I) &=& - \frac{1}{2} S_{\bigtriangleup} \cdot I + 2 \langle J, 
\Ric^-_{\bigtriangleup} \rangle K , \label{gl-57} \\[0.5em]
\tilde{R} (K) &=& - \frac{1}{2} S_{\Box} \cdot K + 2 \langle J, 
\Ric^-_{\bigtriangleup} \rangle I  \ .  \label{gl-58}
\end{eqnarray}

Since $\Ric^-_{\star}, \Ric^-_{\bigtriangleup}$ and $\Ric^-_{\Box}$ are skew 
symmetric and, moreover, it holds that (compare (14))
\begin{equation*} \label{gl-59} 
\{ \Ric^-_{\star} , J\} = \{ \Ric^-_{\bigtriangleup} , I\} = \{ \Ric^-_{\Box} , K\}=0 \ , 
\end{equation*}

analogous considerations using (\ref{gl-45}) yield the equations
\begin{eqnarray*}
\Ric^-_{\star} &=& \langle I,  \Ric^-_{\star} \rangle I + \langle K , 
\Ric^-_{\star} \rangle K , \label{gl-60} \\[0.5em]
\Ric^-_{\bigtriangleup} &=& \langle J,  \Ric^-_{\bigtriangleup} \rangle J - 
\langle K , \Ric^-_{\star} \rangle K , \label{gl-61} \\[0.5em]
\Ric^-_{\Box} &=& - \langle I,  \Ric^-_{\star} \rangle I - \langle J , 
\Ric^-_{\bigtriangleup} \rangle J , \label{gl-62} 
\end{eqnarray*}

and, hence, the local relation
\begin{equation} \label{gl-63} 
| \Ric^-_{\bigtriangleup} |^2 + | \Ric^-_{\Box} |^2 - | \Ric^-_{\star} |^2 = 
2 \langle J, \Ric^-_{\bigtriangleup} \rangle^2 \ . 
\end{equation}

Equation (\ref{gl-20}) shows that $\tilde{R} (A)=0$ for all $A$ with 
$[A,J]=0$. Thus, locally we have
\begin{equation} \label{gl-64} 
| \tilde{R} |^2 = | \tilde{R} (I) |^2 + | \tilde{R} (K) |^2 \ . 
\end{equation}

By (\ref{gl-57}), (\ref{gl-58}), (\ref{gl-64}), we find
\begin{equation} \label{gl-65} 
| \tilde{R} |^2 = \frac{1}{4} (S_{\bigtriangleup}^2 + S_{\Box}^2) +4 (| \Ric^-_{\bigtriangleup} |^2 + | \Ric^-_{\Box} |^2 - | \Ric^-_{\star}|^2 ) \ . 
\end{equation}

By definition of $\tilde{R}^-$, we have generally
\begin{equation*} \label{gl-66} 
\tilde{R}^- (A)= \frac{1}{2} (\tilde{R} (A) - \tilde{R} (JA) \circ J) \ . 
\end{equation*}

Using (\ref{gl-57}), (\ref{gl-58}), we particularly obtain
\begin{eqnarray*}
\tilde{R}^- (I) &=& \frac{1}{4} (S_{\Box} - S_{\bigtriangleup} )I + 2 \langle
J, \Ric^-_{\bigtriangleup} \rangle K , \label{gl-67} \\[0.5em] 
\tilde{R}^- (K) &=& - \frac{1}{4} (S_{\Box} - S_{\bigtriangleup} )K + 2 \langle
J, \Ric^-_{\bigtriangleup} \rangle I , \label{gl-68} \ . 
\end{eqnarray*}

With (\ref{gl-63}) this yields analogously
\begin{equation} \label{gl-69} 
| \tilde{R}^- |^2 = \frac{1}{8} (S_{\bigtriangleup} - S_{\Box} )^2 + 4 
(| \Ric_{\bigtriangleup}^- |^2 + | \Ric^-_{\Box} |^2 - | \Ric^-_{\star} |^2 ) \ . 
\end{equation}

Since the decomposition $\tilde{R} = \tilde{R}^+ + \tilde{R}^-$ is orthogonal,
we have 
\begin{equation} \label{gl-70} 
| \tilde{R} |^2 = | \tilde{R}^+ |^2 + | \tilde{R}^- |^2
\end{equation}

implying the equation
\begin{equation} \label{gl-71} 
| \tilde{R}^+ |^2 = \frac{1}{8} (S_{\star} - S)^2
\end{equation}

by (\ref{gl-46}), (\ref{gl-65}) and (\ref{gl-69}). Using (\ref{gl-46}), (\ref{gl-54})
and (\ref{gl-69}) we immediately obtain the following assertion.

\begin{lem} \label{lem-2-3}
For any almost Hermitian $4$-manifold, the equation
\begin{equation} \label{gl-72} 
| W_+ |^2 = \frac{3}{8} (S_{\star} - \frac{S}{3} )^2 + 8 | \Ric_{\star}^- |^2 + |
\tilde{R}^- |^2
\end{equation}

is valid.
\end{lem}

We denote by $P_1, P_2 : \Lambda^2_+ \to \Lambda^2_+$ the projection corresponding to the splitting (\ref{gl-34}) with $P_1 (\Lambda^2_+)= \R  \cdot \Omega$ and 
$P_2 (\Lambda^2_+)= P_+ (\Omega^{\perp})$. Moreover,  we shall use the 
notations
\begin{displaymath}
\lambda := \frac{1}{4} (S_{\star} - \frac{S}{3} ) \quad , \quad F:= \lambda (
2 P_1- P_2) \quad , \quad G:= W_+ - F \ . 
\end{displaymath}

\begin{lem} \label{lem-2-4}
The equation
\begin{equation} \label{gl-73} 
| W_+ |^2 = 6 \lambda^2 + | G |^2
\end{equation}

holds for any almost Hermitian $4$-manifold $(M,g,J)$. 
\end{lem}

\begin{proof}
Let $(I,K)$ be any local quaternionic supplement of $J$ and let $( \Omega^J, \Omega^I,
\Omega^K)$ be the coframe of the negatively oriented frame $(\Omega_J, \Omega_I, 
\Omega_K)$ of $\Lambda^2_+$. Using the local representations
\begin{equation*} \label{gl-74} 
P_1 = ỏ\Omega^J \otimes \Omega_J \quad , \quad P_2 = \Omega^I \otimes
\Omega_I + \Omega^K \otimes \Omega_K
\end{equation*}

we calculate
\begin{displaymath}
\langle W_+ , P_1 \rangle = \tr (W_+ \circ P_1) = \tr (W_+ \circ (\Omega^J \otimes
\Omega_J))= 
\end{displaymath}
\begin{displaymath}
\tr (\Omega^J \otimes W_+ (\Omega_J))= \Omega^J (W_+ (\Omega_J)) 
\stackrel{\mbox{(\ref{gl-51})}}{=} \frac{1}{2} (S_{\star} - \frac{S}{3}) = 2 \lambda
\end{displaymath}

and, furthermore,
\begin{displaymath}
\langle W_+ , P_2 \rangle = \Omega^I (W_+ (\Omega_I)) + \Omega^K (W_+ (\Omega_K)) 
\stackrel{\mbox{(\ref{gl-52}),(\ref{gl-53})}}{=}
\end{displaymath}
\begin{displaymath}
\frac{1}{2} (S_{\bigtriangleup} - \frac{S}{3} ) + \frac{1}{2} (S_{\Box} - \frac{S}{3})
\stackrel{\mbox{(\ref{gl-46})}}{=} - \frac{1}{2} (S_{\star} - \frac{S}{3} )= - 2
\lambda \ . 
\end{displaymath}

Thus, we obtain
\begin{equation} \label{gl-75} 
\langle W_+ , P_1 \rangle = 2 \lambda = - \langle W_+ , P_2 \rangle \ . 
\end{equation}

Finally, we have
\begin{eqnarray*}
&& | G |^2 = \langle W_+ - F , W_+ - F \rangle = | W_+ |^2 - 2 \langle W_+ , F 
\rangle + | F |^2 =\\[0.5em]
&& | W_+ |^2 - 2 \lambda (2 \langle W_+ , P_1 \rangle - \langle W_+ , P_2 \rangle ) + 
| F |^2 \stackrel{\mbox{(\ref{gl-75})}}{=} \\[0.5em]
&& | W_+ |^2 - 12 \lambda^2 + |F|^2 = |W_+ |^2 - 12 \lambda^2 + 6 \lambda^2 = 
| W_+ |^2 - 6 \lambda^2 \ . 
\end{eqnarray*}
\end{proof}

\begin{prop} \label{prop-2-1}
For any almost Hermitian $4$-manifold, the following assertions are equivalent:\\
(i) The characteristic polynomial of the self-dual part $W_+$ of the Weyl tensor has
the form
\begin{equation*} \label{gl-76} 
\chi (t)= (t - 2 \lambda)(t+ \lambda)^2 \ . 
\end{equation*}

(ii) The equation
\begin{equation} \label{gl-77} 
| W_+ |^2 = 6 \lambda^2
\end{equation}

is satisfied.\medskip

(iii) $W_+$ is given by
\begin{equation} \label{gl-78} 
W_+ = \lambda \cdot (2 P_1 - P_2 ) \ . 
\end{equation}

(iv) The vanishing conditions
\begin{equation*} \label{gl-79} 
\Ric^-_{\star} = 0 \quad , \quad \tilde{R}^- =0
\end{equation*}

are fulfilled.
\end{prop}

\begin{proof}
By Lemma \ref{lem-2-3} and Lemma \ref{lem-2-4}, the equivalence of the assertions 
(ii) - (iv) are obvious. Furthermore, (iii) implies (i) and (i) implies (ii). 
\end{proof}

An immediate consequence of this proposition is the following assertion.

\begin{cor} \label{cor-2-1} 
For any almost Hermitian $4$-manifold, condition (\ref{gl-77}) implies the conformally 
invariant property
\begin{equation} \label{gl-80} 
(\mathrm{det} (W_+))^2 = \frac{1}{54} | W_+ |^6
\end{equation}

of the self-dual part $W_+$ of the Weyl tensor.
\end{cor}

It is well-known that, for any oriented Riemannian 4-manifold $(M,g)$, the 
characteristic polynomial of $W_+$ is given by
\begin{equation*} \label{gl-81} 
\chi (t) = t^3 - \frac{1}{2} |W_+|^2 \cdot t - \mathrm{det} (W_+) \ . 
\end{equation*}

Thus, (\ref{gl-80}) is equivalent to the fact that $W_+ : \Lambda^2_+ \to 
\Lambda^2_+$ has exactly two eigenvalues on the subset $M_+ \subseteq M$ on 
which $W_+$ does not vanish, i.e., at most two eigenvalues at every point of $M$.\\

We finish this section by a proposition summarizing some well-known facts concerning
K\"ahler 4-manifolds (\cite{7}, XII).

\begin{prop} \label{prop-2-2}
Let $(M, g,J)$ be any K\"ahler $4$-manifold with scalar curvature $S$ and self-dual part $W_+$ of its Weyl tensor and let $M_+ \subseteq M$ denote the subset on which
$W_+$ does not vanish. Then the following assertions are valid:\\

(i) The equations
\begin{equation} \label{gl-82} 
| W_+ |^2 = \frac{S^2}{6} \ , 
\end{equation}
\begin{equation} \label{gl-83} 
\mathrm{det} (W_+) = \frac{S^3}{108}
\end{equation}

are satisfied and, hence, also equation (\ref{gl-80}). Moreover, 
the characteristic polynomial of $W_+$ has the form 
\begin{equation} \label{gl-84} 
\chi (t) =( t - \frac{S}{3}) (t + \frac{S}{6})^2 \ . 
\end{equation}

(ii) $W_+$ is given by
\begin{equation} \label{gl-85} 
W_+ = \frac{S}{6} (2 P_1 - P_2)
\end{equation}

and it holds that
\begin{equation} \label{gl-86} 
\nabla W_+ = \frac{1}{6} dS \otimes (2 P_1 - P_2) \ , 
\end{equation}
\begin{equation} \label{gl-87} 
| \nabla W_+ |^2 = \frac{1}{6} | \nabla S |^2 \ . 
\end{equation}

(iii) On $M_+$ we have the equation
\begin{equation} \label{gl-88} 
| \nabla W_+ |^2 = | \nabla | W_+ ||^2 \ . 
\end{equation}
\end{prop}

\begin{proof}
On any K\"ahler manifold, $\Ric_{\star} = \Ric$ and hence
$S_{\star} = S$,  $\Ric^-_{\star} =0$. Moreover, in the K\"ahler case, we have $\tilde{R} =0$. Hence, (45) becomes 
(\ref{gl-82}). (\ref{gl-73}), (\ref{gl-82}) and $S_{\star} =S$  imply $G=0$ and  hence 
equation (\ref{gl-85}). This equation immediately yields (\ref{gl-83})
and (\ref{gl-84}). Moreover, by $P_1 + P_2 = P_+ , \nabla
P_+ =0$ and $\nabla J=0$, we have 
\begin{equation} \label{gl-93} 
\nabla P_1 =0 \quad , \quad \nabla P_2 =0 \ . 
\end{equation}

Using (\ref{gl-93}) we obtain (\ref{gl-86}) from (\ref{gl-85}). (\ref{gl-86})
implies (\ref{gl-87}). Finally, we find (\ref{gl-88}) using (\ref{gl-82}) and
(\ref{gl-87}).
\end{proof}

\section{Integrability conditions for almost Hermitian $4$-manifolds.}
In the following  $M_+$ always denotes the open subset of $M$ on which $W_+$ does not
vanish. Further on we use the notation $\lambda := \frac{1}{4}(S_{\star} - 
\frac{S}{3})$.

\begin{lem} \label{lem-3-1}
For any almost Hermitian $4$-manifold $(M,g,J)$, equation (\ref{gl-77}) implies the equation
\begin{equation*} \label{gl-94} 
| \nabla W_+ |^2 = | \nabla | W_+ ||^2 + \frac{3}{2} | W_+|^2 \cdot | \nabla P_1 |^2 
\end{equation*}

on the subset $M_+ \subseteq M$.
\end{lem}

\begin{proof}
By Proposition \ref{prop-2-1}, (\ref{gl-77}) is equivalent to (\ref{gl-78}). Using
$P_1 + P_2 = P_+$ equation (\ref{gl-78}) becomes
\begin{equation*} \label{gl-95} 
W_+ = \lambda (3 P_1 -P_+) \ . 
\end{equation*}

Since $P_+$ is parallel, applying $\nabla_X$ to this equation we obtain
\begin{equation} \label{gl-96} 
\nabla_X W_+ = X (\lambda) \cdot (3 P_1 - P_+)  + 3 \lambda \nabla_X P_1
\end{equation}

for any vector field $X$. Moreover, on $M_+$, (\ref{gl-77}) yields 
\begin{equation} \label{gl-97} 
| \nabla | W_+ ||^2 = 6 | \nabla \lambda |^2 \ . 
\end{equation}

The equations $\langle P_1, P_1 \rangle = \langle P_+, P_1 \rangle =1, \nabla
P_+ =0$ imply 
\begin{equation} \label{gl-98} 
\langle P_1 , \nabla_X P_1 \rangle = \langle P_+ , \nabla_X P_1 \rangle =0 \ . 
\end{equation}

Finally, we calculate
\begin{displaymath}
| \nabla W_+ |^2 = \langle \nabla_{X_k} W_+ , \nabla_{X^k} W_+ \rangle
= \tr (\nabla_{X_k} W_+ \circ \nabla_{X^k} W_+ ) \stackrel{\mbox{(\ref{gl-96})}}{=}
\end{displaymath}
\begin{displaymath}
6 | \nabla \lambda |^2 + 6 \lambda \langle 3P_1 - P_+ , \nabla_{\nabla \lambda}
P_1 \rangle + 9 \lambda^2 | \nabla P_1 |^2 \stackrel{\mbox{(\ref{gl-77}), (\ref{gl-97}), (\ref{gl-98})}}{=}
\end{displaymath}
\begin{displaymath}
| \nabla | W_+ ||^2 + \frac{3}{2} |W_+ |^2 \cdot | \nabla P_1 |^2 \ . 
\end{displaymath}
\end{proof}

\begin{thm} \label{thm-3-1}
Let $(M,g,J)$ be any almost Hermitian $4$-manifold such that $M_+$ is dense in $M$. Then,
$(M,g,J)$ is K\"ahler if and only if the equation (\ref{gl-77}) is satisfied on $M$ and, moreover, equation (\ref{gl-88}) on $M_+$.
\end{thm}

\begin{proof}
By Lemma \ref{lem-3-1}, (\ref{gl-88}) implies $\nabla P_1 =0$ on $M_+$ and, hence, on 
$M$ since $M_+$ is dense in $M$. On any section $\Omega_A$ of $\Lambda^2_+$, $P_1$ acts
by $P_1 (\Omega_A)= \langle J, A \rangle \Omega_J$. Hence, 
$\nabla P_1 =0$ is equivalent to $\nabla J=0$. Conversely, in the K\"ahler case, Proposition \ref{prop-2-2} shows that we have (\ref{gl-77}) and 
(\ref{gl-88}) then.
\end{proof}

An immediate consequence of this theorem is the following assertion.

\begin{cor} \label{cor-3-1}
An almost Hermitian $4$-manifold $(M,g,J)$ with $\mathrm{supp} (W_+ )=M$ is
K\"ahler if and only if the equations $S_{\star} = S, |W_+|^2 = S^2 /6$ are 
valid and on $M_+$ the equation $|\nabla W_+ | = | \nabla |W_+ ||$.
\end{cor}

We remark that the almost complex structure $J$ enters the integrability
conditions of Theorem \ref{thm-3-1} or Corollary \ref{cor-3-1}, respectively, only
via $S_{\star}$. In particular, these conditions do not contain any derivatives
of $J$. \\

In the following we combine the basic condition (\ref{gl-77}) of this section with the divergence condition
\begin{equation*} \label{gl-99} 
\delta W_+ =0
\end{equation*}
and also with the condition
\begin{equation} \label{gl-100} 
\delta W_+ + \nabla \log | W_+ | \lrcorner W_+ =0
\end{equation}

on $M_+$. The proof of Theorem \ref{thm-3-3} shows that the last equation is satisfied on every K\"ahler 4-manifold.\\

We recall that $\delta W$ is an endomorphism valued $1$-form (tensor field of
type (2,1)) locally defined by $\delta W(X):= (\nabla_{X_k} W)(X,X^k)$. 
Moreover, for any vector field $X \in \Gamma (TM)$, by $X \lrcorner W$ we denote
the endomorphism valued $1$-form  given by $(X \lrcorner W)(Y):= W(X,Y)$. Using the
identification introduced in section 1 we have the relation $W(X,Y) = \frac{1}{2}
W(g(Y) \otimes X-g(X) \otimes Y)$, where the skew symmetric endomorphism $g(Y) \otimes
X - g(X) \otimes Y$ acts by $(g(Y) \otimes X - g(X) \otimes Y)(Z) :=
g(Y,Z)X - g(X, Z)Y$. Locally we obtain
\begin{equation*} \label{gl-101} 
\delta W(X) = \frac{1}{2} (\nabla_{X_k} W)( \xi^k \otimes X - g(X) \otimes X^k)
\ , 
\end{equation*}

where $(\xi^1 , \ldots , \xi^n)$ is the coframe corresponding to the local frame
$(X_1 , \ldots , X_n)$. For any oriented Riemannian $4$-manifold $(M,g)$, the 
decomposition $W=W_+ \oplus W_-$ implies the decomposition
\begin{equation*} \label{gl-102} 
\delta W= \delta W_+ \oplus \delta W_- \ , 
\end{equation*}

where $\delta W_+$ and $\delta W_-$ are locally given by 
\begin{equation} \label{gl-103} 
\delta W_{\pm} (X)= \frac{1}{2} (\nabla_{X_k} W_{\pm} )( \xi^k \otimes
X - g(X) \otimes X^k ) \ . 
\end{equation}

\begin{lem} \label{lem-3-2} 
Let $(M, g,J)$ be an almost Hermitian $\mathrm{4}$-manifold satisfying condition (\ref{gl-77}).
Then, for any vector field $X$ and any local quaternionic supplement $(I,K)$ of $J$, we
have the equation
\begin{equation} \label{gl-104} 
\delta W_+ (X) = - \frac{1}{4} (( 3 \lambda \delta \Omega (X) + 2 d
\lambda (JX))J + 
3 \lambda \cdot \nabla_{JX} J - d \lambda (IX)I - d \lambda (KX)K) \ . 
\end{equation}
\end{lem}

\begin{proof}
We already know that (\ref{gl-77}) implies equation (\ref{gl-96}) being equivalent
to
\begin{equation*}
\nabla_X W_+ =X (\lambda)(2 P_1 - P_2) + 3 
\lambda \nabla_X P_1 \ . 
\end{equation*}

Moreover, for any skew symmetric endomorphisms $A$ and any local quaternionic 
supplement$(I,K)$ of $J$, it holds that 
\begin{displaymath}
P_1 (A) = \langle J,A \rangle J \quad , \quad P_2 (A)= \langle I, A \rangle I 
+ \langle K,A \rangle K \ , 
\end{displaymath}
\begin{displaymath}
(\nabla_X P_1)(A) = \langle \nabla_X J , A \rangle J + \langle J,A \rangle \nabla_X
J \ . 
\end{displaymath}

Hence, by the former expression for $\nabla_X W_+$, we obtain
\begin{equation*}
\begin{array}{l}
(\nabla_X W_+)(A)=(2X (\lambda) \langle J,A \rangle + 3 \lambda \langle
\nabla_X J, A \rangle )J +\\[0.5em]
3 \lambda \langle J,A \rangle \nabla_X J - X(\lambda)(\langle I,A \rangle
I + \langle K,A \rangle K) \ . 
                             \end{array}
\end{equation*}

Now, by a direct calculation using (\ref{gl-103}), we obtain 
(\ref{gl-104}).
\end{proof}

\begin{thm} \label{thm-3-2}
Let $(M,g,J)$ be an almost Hermitian $4$-manifold with the properties (\ref{gl-77})
and $\mathrm{supp} (W_+)=M$. Furthermore, we suppose that at least one of the 
following conditions is satisfied:\\
\mbox{} \quad (i) $(M,g)$ is Einstein.\\
\mbox{} \quad (ii) The Ricci tensor is parallel $(\nabla \Ric =0)$.\\
\mbox{} \quad (iii) The curvature tensor is harmonic $(\delta R=0)$.\\
\mbox{} \quad (iv) The Weyl tensor is harmonic $(\delta W=0)$.\\
\mbox{} \quad (v) $W_+$ is divergence free $(\delta W_+ =0)$.\\[0.4em]
Then the following assertions are valid:
\begin{itemize}
\item[(a)] $J$ is integrable, i.e., $(M,g,J)$ is a Hermitian manifold.
\item[(b)] On $M_+ \subseteq M$, $g$ is conformally equivalent to a 
K\"ahler metric $\bar{g}$ such that $(M_+ , \bar{g}, J)$ is a K\"ahler manifold.
\item[(c)] $(M,g,J)$ itself is a K\"ahler manifold if and only if the function
$| W_+|$ is constant.
\end{itemize}
\end{thm}

\begin{proof}
First of all, we note the implications (i) $\to$ (ii) $\to$ (iii) $\to$ (iv)
$\to$ (v). Thus, it suffices to consider the case of condition (v). By Lemma
\ref{lem-3-2}, on $M_+ \subseteq M$, the equation $\delta W_+ =0$ implies the 
equations
\begin{equation} \label{gl-105} 
\delta \Omega (X) = - \frac{2}{3 \lambda} d \lambda (JX) \ , 
\end{equation}
\begin{equation} \label{gl-106} 
\nabla_{JX} J= \frac{1}{3 \lambda} (d \lambda(IX) I + d \lambda (KX) K) \ . 
\end{equation}

We see that (\ref{gl-106}) implies (\ref{gl-105}) and, moreover, the relation
\begin{equation} \label{gl-107} 
\nabla_{JX} J= J \circ \nabla_X J \ . 
\end{equation}

It is known that (\ref{gl-107}) forces the vanishing of the Nijenhuis tensor $N_J$.
Thus, we have $N_J =0$ on $M_+$ and, hence, on $M$ since $M_+$ is dense in $M$. 
This proves (a). By (\ref{gl-77}), assertion (c) is an immediate consequence
of (\ref{gl-106}). Finally, we prove assertion (b). Let us consider a conformal
transformation $\bar{g} = e^f \cdot g$. Then, the covariant derivatives $\bar{\nabla}$
of $\bar{g}$ and $\nabla$ of $g$ are related by
\begin{equation*} \label{gl-108} 
\bar{\nabla}_X = \nabla_X + \frac{1}{2} (X(f) +df \otimes X - g (X) \otimes
\nabla f) \ . 
\end{equation*}

This yields
\begin{equation} \label{gl-109} 
\bar{\nabla}_X J= \nabla_X J + \frac{1}{2} [df \otimes X - g(X) \otimes \nabla
f , J ] \ . 
\end{equation}

On the other hand, choosing any local quaternionic supplement $(I,K)$ of $J$ we find 
the relation 
\begin{equation} \label{gl-110} 
[df \otimes X - g(X) \otimes \nabla f , J] =df (KX)I - df (IX)K \ . 
\end{equation}

Inserting (\ref{gl-110}) into (\ref{gl-109}) we obtain
\begin{equation} \label{gl-111} 
\bar{\nabla}_X J= \nabla_X J + \frac{1}{2} df (KX) I - \frac{1}{2} df (IX)K\ . 
\end{equation}

Hence, the special choice
\begin{displaymath}
e^f = \sqrt[3]{\lambda^2} \stackrel{\mbox{(\ref{gl-77})}}{=} \sqrt[3]{\frac{1}{6}
| W_+|^2}
\end{displaymath}

yields
\begin{displaymath}
\bar{\nabla}_X J \stackrel{\mbox{(\ref{gl-111})}}{=} \nabla_X J + 
\frac{1}{3 \lambda} (d \lambda (KX)I- d \lambda (IX)K) 
\stackrel{\mbox{(\ref{gl-106})}}{=} 0 \ . 
\end{displaymath}
\end{proof}

\begin{prop} \label{prop-3-1} 
For any almost Hermitian $4$-manifold $(M,g,J)$ which satisfies condition (\ref{gl-77}), 
the equation 
\begin{equation} \label{gl-112} 
\delta W_+ (X)= - \frac{3}{4} \lambda (\delta \Omega (X) J + 
\nabla_{JX} J) - ( \nabla \log | W_+ | \lrcorner W_+ )(X)
\end{equation}

is valid for any vector field $X$ on $M_+$.
\end{prop}

\begin{proof}
Using any local quaternionic supplement $(I,K)$ of $J$ we calculate 
\begin{displaymath}
(\nabla \log | W_+ |  \lrcorner W_+)(X)= W_+ (\nabla \log | W_+ | , X)
\stackrel{\mbox{(\ref{gl-77})}}{=} W_+ (\nabla \log | \lambda | , X)= 
\end{displaymath}
\begin{displaymath}
W_+ (\frac{\nabla \lambda}{\lambda}, X)= \frac{1}{2} W_+ (g (X) \otimes \frac{\nabla \lambda}{\lambda} - \frac{d \lambda}{\lambda} \otimes X) \stackrel{\mbox{(\ref{gl-78})}}{=}
\end{displaymath}
\begin{displaymath}
\frac{1}{2}(2 P_1 - P_2)(g(X) \otimes \nabla \lambda - d \lambda \otimes X)=
\end{displaymath}
\begin{displaymath}
-2 \langle d \lambda \otimes X, J \rangle J + \langle d \lambda \otimes
X,I \rangle I + \langle d \lambda \otimes X, K \rangle K =
\end{displaymath}
\begin{displaymath}
\frac{1}{4} (2 d \lambda (JX)J - d \lambda (IX) - d \lambda (KX)K) \ . 
\end{displaymath}

Thus, we have
\begin{equation} \label{gl-113} 
(\nabla \log | W_+ | \lrcorner W_+)(X)= \frac{1}{4} (2 d \lambda (JX) J - d \lambda (IX) - d \lambda (KX)K) \ . 
\end{equation}

Inserting this into (\ref{gl-104}) we obtain (\ref{gl-112}).
\end{proof}

We remark that the tensor field $\nabla \log | W_+ | \lrcorner W_+$ is defined
not only on $M_+$. It can be extended to M by equation (\ref{gl-113}). 

\begin{thm} \label{thm-3-3}
An almost Hermitian $4$-manifold $(M,g,J)$ with $\mathrm{supp} (W_+)=M$ is K\"ahler
if and only if it satisfies the conditions (\ref{gl-77}) and (\ref{gl-100}).
\end{thm}

\begin{proof}
By Proposition \ref{prop-3-1}, the conditions (\ref{gl-77}), (\ref{gl-100}) immediately
imply $\nabla J=0$. Conversely, by Proposition \ref{prop-2-2}, we have equation
(\ref{gl-86}) in the K\"ahler case. Using (\ref{gl-86}) we obtain the equation
\begin{equation} \label{gl-114} 
\delta W_+ + \nabla \log |S| \lrcorner W_+ =0
\end{equation}

and, hence, equation (\ref{gl-100}) by (\ref{gl-82}).
\end{proof}

Since (\ref{gl-82}) and $S_{\star} = S$ imply (\ref{gl-77}), the following  corollary
is an immediate consequence of Theorem \ref{thm-3-3}.

\begin{cor} \label{cor-3-2}
An almost Hermitian $4$-manifold $(M,g,J)$ with $\mathrm{supp} (S)=M$ is 
K\"ahler if and only if the equations $S_{\star} =S, |W_+ |^2 =S^2/6$ are valid
and equation (\ref{gl-114}) on $M_+$.
\end{cor}

We remark that condition (\ref{gl-114}) is satisfied automatically if the curvature
tensor is harmonic $(\delta R=0)$ and, hence, if the Ricci tensor is parallel or if
the manifold is Einstein.\\

We finish this section by a theorem which is closely related to a result of 
Derdzinski (\cite{6}, 16.67) and which shows that the local K\"ahler property of an 
oriented Riemannian 4-manifold $(M,g)$ with $W_+ \not= 0$ everywhere can be 
formulated without using any almost complex structure. 

\begin{thm} \label{thm-3-4} 
Let $(M,g)$ be an oriented Riemannian $4$-manifold such that $W_+$ vanishes nowhere.
Then the following conditions are equivalent:\\
(i) $g$ is locally a K\"ahler metric with respect to local complex structures that 
are compatible with the given orientation.\\
(ii) The condition (\ref{gl-80}) and (\ref{gl-88}) are satisfied.\\
(iii) The equations (\ref{gl-80}) and (\ref{gl-100}) are valid.
\end{thm}

\begin{proof}
The assumptions $M_+ =M$ and (\ref{gl-80}) imply that $W_+$ has exactly two eigenvalues at every point of $M$. This yields the orthogonal splitting
\begin{equation} \label{gl-115} 
\Lambda^2_+ = E_1 \oplus E_2
\end{equation}

of $\Lambda^2_+$ into the corresponding eigenbundles $E_1, E_2$ of $W_+$ 
with $\mathrm{rank} (E_{\alpha})= \alpha \ (\alpha =1,2)$. By assertion (ii) of Lemma
\ref{lem-2-1},  the local unit sections of $E_1$ define local almost complex
structures compatible with the metric $g$ and the given orientation. Thus, any oriented
Riemannian 4-manifold $(M,g)$ with $M_+ =M$ and (\ref{gl-80}) is locally almost 
Hermitian. Moreover, by $\tr (W_+)=0$, (\ref{gl-115}) implies that $W_+ = \mu (2
Q_1 -Q_2)$, where $Q_1 , Q_2 : \Lambda^2_+ \to \Lambda^2_+$ denote the corresponding projections $(Q_{\alpha} (\Lambda^2_+)= E_{\alpha}, \alpha =1,2)$. Now, the proofs
of Theorem \ref{thm-3-1} and Theorem \ref{thm-3-3} show that, in this case, condition 
(\ref{gl-88}) and also condition (\ref{gl-100}) forces the K\"ahler property of all these local almost Hermitian structures. The converse is true by Proposition 
\ref{prop-2-2} and Theorem \ref{thm-3-3}.
\end{proof}

\section{Integrability conditions for almost K\"ahler $4$-manifolds.}

By (\ref{gl-32}) and (\ref{gl-70})-(\ref{gl-72}), we immediately obtain the 
following assertion.

\begin{prop} \label{prop-4-1} 
On every almost K\"ahler $4$-manifold $(M,g,J)$, the equation
\begin{equation} \label{gl-116} 
| W_+ |^2 - \frac{S^2}{6} = S | \nabla J |^2 + | \nabla J |^4 + 8 | \Ric^-_{\star}
|^2 + | \tilde{R} |^2
\end{equation}

is valid.
\end{prop}

Together with assertion (i) of Proposition \ref{prop-2-2} this proves the 
following result.

\begin{thm} \label{thm-4-1} 
Let $(M,g,J)$ be any almost K\"ahler $4$-manifold with scalar curvature $S \ge 0$.
Then $(M,g,J)$ is K\"ahler if and only if the equation $| W_+ |^2 =S^2 /6$ is
satisfied.
\end{thm}

\begin{cor} \label{cor-4-1} 
Every almost K\"ahler $4$-manifold with the property $S= | W_+ | \sqrt{6}$ is already a K\"ahler manifold.
\end{cor}

For any compact almost Hermitian manifold $(M,g,J)$ we consider the number
\begin{displaymath}
Q(J) := \int_M q (J) \omega \ , 
\end{displaymath}

where $\omega$ denotes the volume form and $q (J)$ the function locally defined
by 
\begin{displaymath}
q(J) := g (( \nabla^2_{X_k, X_l} \Ric )JX^k , JX^l ) \ . 
\end{displaymath}

$Q(J)$ is an obstruction to the K\"ahler property for any compact almost Hermitian manifold (\cite{12}, Section 2). The following
proposition is a corollary of Proposition 2.5 in \cite{12}.

\begin{prop} \label{prop-4-2} 
For any compact almost K\"ahler $4$-manifold $(M,g,J)$, we have the equation
\begin{equation} \label{gl-117} 
Q(J) + \int_M (| \tilde{R}|^2 + 4 | \Ric_{\star}^- |^2 + 
\frac{S}{2} | \nabla J |^2 + | \nabla J |^4 ) \omega =0 \ . 
\end{equation}
\end{prop}

We remark that, in contrast to this paper, in \cite{12} the length of an endomorphism
$A$ is defined by $|A|^2 := \tr (A^* \circ A)$. This explains the factor 4
in equation (\ref{gl-117}). We also have to take into account the different 
definitions of the function $| \nabla J |^2$ in \cite{12} and in this paper.\\
Combining (\ref{gl-116}) and (\ref{gl-117}) we obtain the next proposition.

\begin{prop} \label{prop-4-3} 
On any compact almost K\"ahler $4$-manifold $(M,g,J)$, the equation
\begin{equation} \label{gl-118} 
Q(J) + \frac{1}{2} \int_M (| \tilde{R}|^2 + | \nabla J |^4 + | W_+ |^2 - 
\frac{S^2}{6}) \omega =0
\end{equation}

holds.
\end{prop}

The Weitzenb\"ock formula (\ref{gl-118}) immediately yields the following 
result.

\begin{thm} \label{thm-4-2}
A compact almost K\"ahler $4$-manifold $(M,g,J)$ is K\"ahler if and only if the 
equations $Q(J)=0$ and $| W_+ |^2 =S^2/6$ are satisfied.
\end{thm}

By Remark 3.9 in \cite{12}, we know several conditions that imply the vanishing
 of $Q(J)$. For example, the assumption $[\nabla \Ric, J]=0$, i.e., 
$[\nabla_X \Ric, J]=0$ for all vector fields $X$, implies $Q(J)=0$ and also the 
supposition $[\nabla^2 \Ric, J] =0$, i.e., $[ \nabla^2_{X,Y} \Ric , J]=0$
for all vector fields $X,Y$. In particular, we have $Q(J)=0$ if the Ricci tensor
is parallel $(\nabla \Ric =0)$. Another condition that forces $Q(J)=0$ is
$\delta W=0$ (harmonic Weyl tensor). Thus, the next corollary is an immediate 
consequence of Theorem \ref{thm-4-2}.

\begin{cor} \label{cor-4-2}
Let $(M, g,J)$ be a compact almost K\"ahler $4$-manifold satisfying at least one 
of the following conditions:\\
(i) The relation $[\nabla^2 \Ric , J]=0$ is valid.\\
(ii) It holds that $[\nabla \Ric , J]=0$.\\
(iii) The Weyl tensor is harmonic.\\
(iv) The Ricci tensor is parallel.\\
(v) $(M,g)$ is Einstein.\\
Then $(M,g,J)$ is a K\"ahler manifold if and only if $| W_+ |^2 =S^2 /6$. 
\end{cor}

It is known that the Ricci tensor of a K\"ahler manifold with harmonic Weyl tensor
is parallel (\cite{6}, 16.30). This shows that, in the context of
Corollary \ref{cor-4-2}, the conditions (iii) and (iv) are equivalent. In the 
following we prove a result which is an essential generalization of this part
of Corollary \ref{cor-4-2}.\\
For any Riemannian $n$-manifold $(M,g)$ with $n \ge 3$, the relation
\begin{equation*} \label{gl-119} 
g(\delta W(X) Y, Z)= \frac{n-3}{n} g(X, (\nabla_Y \Ric)Z - (\nabla_Z \Ric)Y
- \frac{1}{2(n-1)} (Y(S)Z-Z(S)Y))
\end{equation*}

is well-known. Using this we obtain the equation
\begin{equation} \label{gl-120} 
\begin{array}{c}
\langle \delta W(X), A \rangle = \frac{n-3}{n(n-2)} g(X, (\nabla_{X_k} \Ric
)AX^k - (\nabla_{AX_k} \Ric)X^k \\[0.5em]
- \frac{1}{2(n-1)} (A(\nabla S) - dS (AX_k)X^k))
\end{array}
\end{equation}

for any vector field $X$ and any skew symmetric endomorphism $A:TM \to TM$.

\begin{lem} \label{lem-4-1} 
Let $(M,g)$ be any oriented Riemannian $4$-manifold. Then, for any vector field 
$X$ and any positively oriented local orthonormal frame $(\Omega_I, \Omega_J,
 \Omega_K)$ of $\Lambda^2_+$, we have the equation
\begin{equation} \label{gl-121} 
\begin{array}{c}
\delta W_+ (X) = \frac{1}{4} (g(X, (\nabla_{X_k} \Ric)IX^k )I +g(X, (\nabla_{X_k}
\Ric )JX^k)J +g (X, (\nabla_{X_k} \Ric)KX^k)K) \\[0.5em]
+ \frac{1}{24} (dS(IX)I+dS (JX) J + dS(KX)K) \ . 
\end{array}
\end{equation}
\end{lem}

\begin{proof}
By definition, it holds that 
\begin{displaymath}
\delta W_+ (X)= P_+ (\delta W(X)) \ . 
\end{displaymath}

Thus, we obtain the local representation
\begin{displaymath}
\delta W_+ (X)= \langle \delta W(X), I \rangle I+ \langle \delta W(X), J \rangle
J+ \langle \delta W(X), K \rangle K \ , 
\end{displaymath}

where $(\Omega_I, \Omega_J, \Omega_K)$ is any positively oriented local 
orthonormal frame of $\Lambda^2_+$. By (\ref{gl-120}), this yields
\begin{eqnarray*}
\delta W_+ (X) &=& \frac{1}{4} g (X, (\nabla_{X_k} \Ric) IX^k - 
\frac{1}{6} I(\nabla S)) I\\[0.5em]
&& + \frac{1}{4} g(X, \nabla_{X_k} \Ric)JX^k - \frac{1}{6} J(\nabla S))J
\\[0.5em]
&& + \frac{1}{4} g (X, (\nabla_{X_k} \Ric)KX^k - \frac{1}{6} K(\nabla S))K
\end{eqnarray*}

and, hence, (\ref{gl-121}).
\end{proof}

Now, for any almost Hermitian 4-manifold $(M,g,J)$, we consider the 
endomorphism valued 1-form $\Theta$ defined by
\begin{displaymath}
\Theta := \frac{1}{6} \nabla S \lrcorner (2 P_1 - P_2) \ , 
\end{displaymath}

where $P_1, P_2$ as before denote the projections of splitting (\ref{gl-34}).
Then, for any vector field $X$, we have locally 
\begin{equation} \label{gl-122} 
\Theta (X)= \frac{1}{24} (2dS(JX)J - dS (IX)I - dS (KX)K) \ , 
\end{equation}

where $(I,K)$ is any local quaternionic supplement of $J$. 

\begin{lem} \label{lem-4-2} 
Let $(M,g,J)$ be any almost Hermitian $4$-manifold. Then, the following assertions
are valid: The condition
\begin{equation*} \label{gl-123} 
\delta W_+ + \Theta =0 \quad \quad (\delta W_+=0)
\end{equation*}

is satisfied if and only if we have the equation
\begin{equation}  \label{gl-124} 
(\nabla_{X_k} \Ric )JX^k = \frac{1}{2} J(\nabla S) \quad ((\nabla_{X_k} \Ric)JX^k = \frac{1}{6} J(\nabla S))
\end{equation}

and, moreover, for any local quaternionic supplement $(I,K)$ of $J$, the 
equations
\begin{equation*} \label{gl-125} 
(\nabla_{X_k} \Ric)IX^k=0 \quad , \quad (\nabla_{X_k} \Ric ) KX^k =0
\end{equation*}
\begin{displaymath}
((\nabla_{X_k} \Ric)IX^k = \frac{1}{6} I(\nabla S) \quad , \quad (\nabla_{X_k}\Ric)KX^k = \frac{1}{6} K(\nabla S)) \ . 
\end{displaymath}
\end{lem}

\begin{proof}
We apply Lemma \ref{lem-4-1} to such local orthonormal frames $(\Omega_I, 
\Omega_J, \Omega_K)$ of $\Lambda^2_+$ for which $(I,K)$ is any local quaternionic 
supplement of $J$. Then, by (\ref{gl-121}) and (\ref{gl-122}), we immediately
obtain the assertion of our lemma.
\end{proof}

Now, for any almost Hermitian manifold $(M,g,J)$, we introduce the vector valued
2-forms $\varphi$ and $\psi$ defined by 
\begin{displaymath}
\varphi (X,Y):= (\nabla_X J) Y - (\nabla_Y J)X \ , 
\end{displaymath}
\begin{displaymath}
\psi (X,Y):= \frac{1}{8} ([\nabla_X \Ric, J]Y- [\nabla_Y \Ric , J]X - 
[\nabla_{JX} \Ric, J] JY + [\nabla_{JY} \Ric , J]JX) \ . 
\end{displaymath}

Using the notation $\langle \varphi , \psi \rangle := \frac{1}{2} g( \varphi
(X^k, X^l), \psi (X_k, X_l))$, by (63) and (72) in \cite{12}, 
we have the following lemma.

\begin{lem} \label{lem-4-3} 
Let $(M, g,J)$ be any almost K\"ahler manifold. Then the equation
\begin{equation} \label{gl-126} 
\langle \varphi , \psi \rangle = \frac{1}{4} g(J(\varphi (X^k , X^l)), (
\nabla_{X_k} \Ric) X_l - (\nabla_{X_l} \Ric)X_k)
\end{equation}

is valid. Moreover, if $M$ is compact, then
\begin{equation*} \label{gl-127} 
Q(J) = 2 \int_M \langle \varphi , \psi \rangle \omega \ . 
\end{equation*}
\end{lem}
\mbox{} \\

If $(M,g,J)$ is any almost Hermitian 4-manifold and $(I,K)$ any quaternionic 
supplement of $J$ on a sufficiently small open subset $U \subseteq M$, then
there are uniquely determined vector fields $\xi, \eta$ on $U$ such that 
\begin{equation} \label{gl-128} 
\nabla_X J= g(\xi, X)I+g( \eta , X)K
\end{equation}

for any vector field $X$ on $U$. It is well-known that we have $d \Omega_J=0$
then if and only if
\begin{equation} \label{gl-129} 
\eta = J \xi \ . 
\end{equation}

\begin{lem} \label{lem-4-4} 
Let $(M,g,J)$ be any almost K\"ahler $\mathrm{4}$-manifold. Then the function $\langle \varphi, 
\psi \rangle$ is locally given by
\begin{equation} \label{gl-130} 
\langle \varphi , \psi \rangle = \frac{1}{2} g((\nabla_{X_k} \Ric ) KX^k , \xi)
- \frac{1}{2} g (( \nabla_{X_k} \Ric)IX^k , \eta) \ , 
\end{equation}

where $(I,K)$ is any local quaternionic supplement of $J$ and $\xi,\eta$ the 
corresponding vector fields defined by (\ref{gl-128}).
\end{lem}

\begin{proof}
Using (82) and (83), we compute
\begin{eqnarray*}
\langle \varphi , \psi \rangle & {=} & \frac{1}{2}
g(J((\nabla_{X^k} J)X^l) , (\nabla_{X_k} \Ric) X_l - ( \nabla_{X_l} \Ric)X_k)
{=} \\[0.5em]
&& \frac{1}{2} g (-g (\xi, X^k)KX^l + g(\eta , X^k) IX^l, (\nabla_{X_k} \Ric)
X_l - ( \nabla_{X_l} \Ric)X_k)=\\[0.5em]
&& - 2 \langle K, \nabla_{\xi} \Ric \rangle + \frac{1}{2} g ((\nabla_{X_l} \Ric)KX^l , \xi)+ 2 \langle I, \nabla_{\eta} \Ric \rangle - \frac{1}{2} g(( \nabla_{X_l} \Ric)IX^l , \eta ) \ . 
\end{eqnarray*}

This yields (\ref{gl-130}) since $\langle K, \nabla_{\xi} \Ric \rangle = \langle I, 
\nabla_{\eta} \Ric \rangle =0$.
\end{proof}

Now, we are able to prove the main results of this section.

\begin{thm} \label{thm-4-3} 
Let $(M,g,J)$ be a compact almost K\"ahler $4$-manifold. Then we have the following:\\
(i) $(M,g,J)$ is K\"ahler if and only if $|W_+ |^2 = S^2/6$ and $\delta W_+ + 
\Theta =0$.\\
(ii) $(M,g,J)$ is a  K\"ahler manifold of constant scalar curvature if and only if
$| W_+ |^2 =S^2/6$ and $\delta W_+ =0$.
\end{thm}

\begin{proof}
By Lemma \ref{lem-4-2} and Lemma \ref{lem-4-4}, we see that $\delta W_+ + \Theta =0$
implies $\langle \varphi , \psi \rangle =0$ and, hence, $Q(J)=0$ by Lemma \ref{lem-4-3}.
Thus, by Theorem \ref{thm-4-2}, $(M,g,J)$ is K\"ahler. Conversely, by Proposition
\ref{prop-2-2} (equations (\ref{gl-82}) and (\ref{gl-86})), we see that the equations
$| W |^2 =S^2/6$ and $\delta W_+ + \Theta =0$ are valid for any K\"ahler manifold.
This proves assertion (i). By Lemma \ref{lem-4-2}, the assumption $\delta W_+ =0$
implies locally
\begin{displaymath}
(\nabla_{X_k} \Ric) IX^k = \frac{1}{6} I(\nabla S) \quad , \quad (\nabla_{X_k} \Ric)
KX^k = \frac{1}{6} K(\nabla S) \ . 
\end{displaymath}

Inserting this into (\ref{gl-130}) we obtain
\begin{displaymath}
\langle \varphi , \psi \rangle = \frac{1}{12} g (K (\nabla S) , \xi) -
\frac{1}{12} g(I(\nabla S) , \eta ) \stackrel{\mbox{(\ref{gl-129})}}{=} 0
\end{displaymath}

and, hence, $Q(J)=0$ as before. Thus, the conditions $\delta W_+ =0$ and $| W_+ |^2 =S^2/6$ also imply the K\"ahler property. But, in the K\"ahler case, equation 
(\ref{gl-124}) in parentheses immediately yields $\nabla S=0$.
\end{proof}

Since the condition $\delta W_+ + \Theta =0$ and also $\delta W_+ =0$ generalizes
the Einstein condition, Theorem \ref{thm-4-3} is related to the Goldberg conjecture 
(see introduction). The following theorems are also results in this direction.

\begin{thm} \label{thm-4-4} 
A compact almost K\"ahler Einstein $4$-manifold is K\"ahler if and only if the 
length of the self-dual part of its Weyl tensor is constant.
\end{thm}

\begin{proof}
By (\ref{gl-72}), we have the inequality
\begin{equation} \label{gl-131} 
| W_+ |^2 \ge \frac{3}{8} (S_{\star} -\frac{S}{3})^2
\end{equation}

for any almost Hermitian 4-manifold. Now, for any compact almost K\"ahler Einstein
4-manifold, it was proved by J. Armstrong \cite{3} that there is at least one 
point at which $S_{\star} =S$. By (\ref{gl-131}), at such a point we have  
\begin{displaymath}
|W_+ |^2 \ge \frac{S^2}{6} \ .
\end{displaymath}

Thus, if $|W_+|$ is constant, this inequality  is satisfied everywhere since $S$ is constant
owing to the Einstein condition. By Proposition \ref{prop-4-3}, this forces $\nabla J=0$ since we have $Q(J)=0$ in the Einstein case. The converse is true by assertion (i) of Proposition \ref{prop-2-2}.
\end{proof}

A theorem of Oguro and Sekigawa \cite{16} asserts that an almost K\"ahler 
Einstein 4-manifold is K\"ahler if and only if its star scalar curvature is constant. 
Theorem \ref{thm-4-4} is a similar result in the compact case.

\begin{cor} \label{cor-4-3}
A compact almost K\"ahler Einstein $4$-manifold is K\"ahler if and only if it has the
property
\begin{equation} \label{gl-132} 
\mathrm{det} (W_+) = \frac{S}{18} | W_+ |^2 \ . 
\end{equation}
\end{cor}

\begin{proof}
By assertion (i) of Proposition \ref{prop-2-2}, we see that equation (\ref{gl-132})
is satisfied in the K\"ahler case. Conversely, the Einstein condition implies $\delta
W_+ =0$ and, hence, by 16.73 in \cite{6}, the equation
\begin{equation} \label{gl-133}
2 | \nabla W_+ |^2 + \bigtriangleup | W_+ |^2 = 18 \mathrm{det} (W_+) - S|W_+|^2 \ . 
\end{equation}

Integrating this equation the supposition (\ref{gl-132}) immediately forces
$\nabla W_+ =0$. In particular, $|W_+|$ is constant. By Theorem \ref{thm-4-4}, this
implies the K\"ahler property.
\end{proof}

\begin{thm} \label{thm-4-5}
A compact almost K\"ahler $4$-manifold with constant negative scalar curvature $S$
is K\"ahler if and only if it satisfies the curvature conditions $\mathrm{det} (W_+) =
S^3 /108$ and $\delta W_+=0$.
\end{thm}

\begin{proof}
By Proposition \ref{prop-2-2}, the conditions of our theorem are satisfied for any 
K\"ahler 4-manifold of constant scalar curvature. Conversely, since $\delta W_+ =0$
and $S$ is constant by supposition, we have
\begin{displaymath}
0 \le 2 \int_M | \nabla W_+ |^2 \omega = \int_M (2 | \nabla W_+ |^2 + \bigtriangleup
| W_+ |^2) \omega \stackrel{\mbox{(\ref{gl-133})}}{=}
\end{displaymath}
\begin{displaymath}
\int_M (18 \mathrm{det} (W_+) - S| W_+ |^2 ) \omega {=}
(-S) \cdot \int_M (|W_+ |^2 - \frac{S^2}{6}) \omega \ . 
\end{displaymath}

Thus, it follows $\int_M (|W_+|^2 - \frac{S^2}{6}) \omega \ge 0$ since $S < 0$. 
By Proposition \ref{prop-4-3}, this yields $\nabla J=0$ since we already
know that $\delta W_+ =0$ implies $Q(J)=0$ in the almost K\"ahler case.
\end{proof}

\begin{cor} \label{cor-4-4}
Let $(M,g,J)$ be a compact almost K\"ahler $4$-manifold with $S <0$ such that at least one of the following conditions is satisfied: \\
(i) The curvature tensor is harmonic $(\delta R=0)$.\\
(ii) The Ricci tensor is parallel.\\
(iii) $(M,g)$ is Einstein.\\
Then $(M,g,J)$ is K\"ahler if and only if $\mathrm{det} (W_+)=S^3 /108$.
\end{cor}

\begin{proof}
We have the implications (iii) $\to$ (ii) $\to$ (i). But $\delta R=0$ implies that
$S$ is constant and $\delta W_+ =0$.
\end{proof}

\vspace{1cm}


\end{document}